\documentclass{amsart}
\usepackage{amsmath}
\usepackage{amssymb}
\usepackage{amscd}
\usepackage[mathscr]{eucal}
\newtheorem{thm}{Theorem}[section]

\theoremstyle{definition}
\newtheorem{defn}[thm]{Definition}
\theoremstyle{remark}

\newtheorem{exa}[thm]{Example}
\numberwithin{equation}{section}

\def\vep{\varepsilon}
\def\tvep{\widetilde{\varepsilon}}

\def\Ki{\widehat{T}_i}
\begin{document}
\title {Asymptotic dynamics of generalized \\ Kantorovich operators }
\author {Krzysztof Bartoszek$^1$ and Wojciech Bartoszek $^2$ }
\address {$^1$Department of Computer and Information Science, Link{\"o}ping University, 
Link{\"o}ping 581 83 Sweden}
\email{$^1$krzysztof.bartoszek@liu.se, krzbar@protonmail.ch}

\email{$^2$BartoszekMathLab@gmail.com}


\date {27 October 2023}
\begin{abstract}
We characterize
the family of 
continuous functions $f\in C([0,1])$ such that
the iterates $\Ki^kf$ converge uniformly on $[0,1]$, where $\Ki$ is a generalized Kantorovich operator. This
gives an affirmative answer  to the problem raised in  2021 by  Acu and Rasa.
\end{abstract}
\maketitle
\pagestyle{myheadings}

\markboth{Bartoszek \& Bartoszek }{Generalized Kantorovich operators}

\section{Introduction}
Let $(X, \varrho )$ be a compact (metric) space and $C(X)$ be the Banach lattice of all real valued continuous
functions $f : X \mapsto \mathbb{R}$,  endowed with the sup-norm $\| f \| = \sup\limits_{x\in X} |f(x) |$. A linear
(continuous) operator $T : C(X) \mapsto C(X)$ is called Markov if $T{\bf 1} = {\bf 1}$ and $Tf \geq 0$ for all
nonnegative $f\in C(X)$. We shall use the
``bracket" dual notation, i.e., $ Tf(x) = T'\delta_x (f) = \langle f , T'\delta_x\rangle \ , $
for all $f\in C(X)$ and all $x\in X$. $T'$ is the dual operator acting on the Banach lattice $C'(X) = M(X)$ of all
$\sigma-$additive finite (Borel) measures on $X$. Clearly $\delta_x $ stands for the Dirac measure concentrated at $x\in X$.
The formula $ Tf(x) = T'\delta_x (f)$ naturally extends $T$ to Borel (bounded) functions $f : X \mapsto \mathbb{R}$.
We underline that $X \ni x \mapsto T'\delta_x \in \mathcal{P}(X)$ is $\varrho-*weak$ continuous
($T$ is Feller and in the case when the phase space $X$ is locally compact the continuity assumption is understood as
weak measure convergence continuity).
The topic of Markov operators (introduced above) is quite well studied
and from the perspective of our work
may be attributed to R. Foguel ([7], [8]), B. Jamison ([12]), R. Sine ([20], [22]),  M. Lin ([16], [17]),  A-M Acu and I. Rasa ([1], [2]) and
[3], [5].  We add that a Markov operator $T : C(X) \mapsto C(X)$ is a positive linear contraction and its operator norm
$\| T \| = 1$.
Markov operators constitute an alternative  functional approach to the notion of stochastic Markov processes,
in particular homogeneous ones.
The  main question in the theory of Markov operators is to characterize the
asymptotic behaviour of the iterates $T^kf$, $k = 0, 1, \dots $,
where $f$ is a fixed function. Turning to stochasticity  let us remark that $T^kf(x)$ is the mean value of a corresponding
Markov process, at the
time instant $t= k$, for an observable $f$ (continuous or measurable).
More precisely, if  \ $\xi_0, \xi_1, \dots $ \
is the associated Markov process, then $T^kf(x) = \mathbb{E}(f(\xi_k ) | \xi_0 = x )$.

We recall that $T$ is said to be {\it strongly ergodic} ({\it mean asymptotically stable}) if for all $f \in C(X)$ the Cesaro averages
$$
A_mf = \frac{f + Tf + \dots T^{m-1}f}{m}
$$
converge uniformly (i.e. in the sup-norm $\| \cdot \|$) to $Pf\in C(X)$. Clearly $P: C(X) \mapsto C(X)$ is a 
Markovian projection. Moreover,
$TP = PT = P$. Hence, $P : C(X) \mapsto C_{{\rm inv}}(X)$, where $C_{{\rm inv}}(X) = \{ f\in C(X) : Tf = f \} $
is the manifold of all $T-$invariant continuous functions. By definition ${\bf 1} \in C_{{\rm inv}}(X) $.
A signed (Radon) measure $\mu \in C(X)' = M(X)$ is called $T-$invariant if $T'\mu = \mu$,
$${\rm i.e. \ \ } \ \ \  \langle f , \mu \rangle = \langle Tf , \mu \rangle = \int_X Tf d\mu =  \int_X f dT'\mu = \langle f , T'\mu \rangle  $$
hold for all $f\in C(X)$. The manifold of all $T-$invariant measures is denoted as $M_{{\rm inv}}(X)$ and the convex *weak-compact
(nonempty) set of $T-$invariant probability measures is denoted as $P_T(X) = \mathcal{P}(X)\cap M_{{\rm inv}}(X)$, where $\mathcal{P}(X)$ is the *weak compact set of all probability measures on $X$ (cf [20]).

It appears that a Markov operator $T$ is strongly ergodic if and only if
$T-$invariant functions separate $T-$invariant measures, i.e.,
for every nonzero $\mu \in M_{{\rm inv}}(X)$ there exists $f\in C_{{\rm inv}}(X)$ such that $\int_X f d\mu = \langle f , \mu \rangle \neq 0$ \
(see [20], [21]). It is worth emphasizing that $T$ is strongly ergodic if and only if
for any pair of distinct $T-$invariant probabilities $\nu \neq \lambda $ there exists an invariant function $f\in C_{{\rm inv}}(X)$ such that
$\int_X f d\nu \neq \int_X fd\lambda$. It is a commonly known fact that  $T-$invariant measures always separate $T-$invariant functions. This is even true  in a much broader context, for all linear contractions on Banach spaces (see [21]). In particular,  for every pair
of invariant functions $f_1\neq f_2 \in C_{{\rm inv}}(X)$ there always exists an invariant probability measure $\nu \in P_T(X)$ such that
$\langle f_1 , \nu \rangle \neq \langle f_2 , \nu \rangle $.

\begin{exa}
Consider the metric space
$X = [0,1]$ with the  standard metric $\varrho(x,y) = |x - y|$. Define $Tf(x) = f(0)(1-x) + f(1)x$. Clearly
$T$ is a Markov projection ($T^2 = T = P$), $C_{{\rm inv}}([0,1]) = {\rm Aff} ([0,1])$  coincides with affine functions on $[0,1]$ and
$M_{{\rm inv}}([0,1]) = \{ \alpha \delta_0 + \beta \delta_1 : \alpha, \beta \in \mathbb{R} \}$. In particular,
$P_T([0,1]) = \{ (1-t)\delta_0 + t\delta_1 :  t \in [0, 1]\} $.
As a projection, this Markov operator is not only strongly ergodic, but {\it uniformly ergodic} (i.e., Cesaro means $A_m $ converge
in the operator norm), or even more, the iterates $T^m = A_m = T$ "converge" in the operator norm to the limit projection.
\end{exa}

The next example is not so trivial as above, and brings us closer
to the class which is a predecessor of our main object of interests.

\begin{exa}
As before, let $X = [0,1]$ and for a fixed natural $k$ we define
$$
B_kf(x) = \sum_{j=0}^k \binom{k}{j} x^j(1-x)^{k-j}f(\tfrac{j}{k}) \ .
$$
The operators
$B_k$ are called standard {\it Bernstein operators}. It is well known that $B_k$ are finite dimensional Markov operators
on $C([0,1])$, hence compact. Asymptotic properties of their iterates $B_k^m$ are well known and it may be proved that
$$
\lim_{m\to \infty} B_k^m f = f(0)(1 - {\bf e}_1) + f(1){\bf e}_1 \ ,
$$
where ${\bf e}_1(x) = x$, and the convergence holds for the operator norm. As $B_k$ is a compact Markov operator, it follows that  (cf [3], [4]) the rate of convergence is geometric (however, this fact can be established using elementary methods only, cf [14], [13], [19],
[9], [10]). As before, the limit projection is $2-$dimensional,  \
$C_{{\rm inv}}([0,1]) = {\rm Aff} ([0,1]) = {\rm lin} \{ {\bf 1}, {\bf e}_1 \}$ \ and \ $M_{{\rm inv}}([0,1]) = \{ \alpha \delta_0 + \beta \delta_1 : \alpha, \beta \in \mathbb{R} \}$.
We notice that
$$
B_k'\delta_x = \sum_{j=0}^k \binom{k}{j}x^j(1-x)^{k-j}\delta_{\tfrac{j}{k}} \ ,
$$
$B_k$ is $k+1$ dimensional and
$\lim\limits_{m\to \infty} B'^m_k \delta_x = (1 - x)\delta_0 + x\delta_1$. Clearly the convergence is in the total variation norm, uniformly for
$x$ and it holds with the geometric rate.
\end{exa}

The following example  is a bit more advanced (it is inspired by R. Sine's ideas, cf [22], Example 1). We recall that a (closed) subset $A$ of a phase space $X$ is called {\it T-invariant} if $T'\delta_a(A) = 1$ for all $a\in A$. By the topological support of a measure $\nu \in \mathcal{P}(X)$ we mean the smallest closed subset $F\subseteq X$ such that $\nu (F) = 1$, or equivalently $\nu (F^c) = 0$ ($|\nu |(F^c) = 0 $ if we
consider a signed measure $\nu \in M(X)$). In this paper measures $\nu \in M(X)$ have  supports, due to the fact that considered here metric spaces are
separable and complete (in the case of general compact Hausdorff spaces $X$ we have a similar situation, once we discuss only Radon measures).

\begin{exa}
Let $X = D_1 = \{ z\in \mathbb{C} : |z| \leq 1 \}$ be the closed unit disc of the complex plain $\mathbb{C}$. Define a Markov operator
$T : C(D_1) \mapsto C(D_1)$ by its transition probabilities (densities) $T'\delta_u $ for \ $u \in D_1$: $T'\delta_u = \delta_u $ \ if \ $ |u| = 1$,
$$
T'\delta_u = \frac{1}{\pi r_u^2} {\bf 1}_{\{ z\in D_1 : | z -  (2 - \frac{1}{|u|})u | \leq r_u \}}, \ \ {\rm if} \ \ \frac{1}{2} \leq |u| < 1 ,
$$
where $r_u = 1 - |u|$, and
$$
T'\delta_u =  \frac{1}{\pi r_u^2} {\bf 1}_{\{  z\in D_1 : |z| \leq r_u\}}, \ \ {\rm if} \ \ 0 \leq |u| \leq \frac{1}{2} \ .
$$
It is easily  seen  that $\delta_u \in P_T(D_1)$ for all $|u| = 1$.
Without detailed explanations we note that the Markov process $\xi_0, \xi_1 , \dots   $ generated by $T$ is irreducible on the interior $D = \{ z\in D_1 : |z| < 1 \}$, i.e.,
for every $z\in D$ and any open subset $U \subseteq D$ there exists a natural $k$ such that
$$
T'^k\delta_z(U)
= T^k{\bf 1}_U(z) = {\rm Prob}(\xi_k \in U | \xi_0 = z) > 0 .
$$
The topological support of $T'\delta_0 $ coincides with $D_1$.
It follows that there exists
another $T-$invariant ergodic probability $\lambda $ concentrated on $D_1$ (stochastically it lives on $D$). By irreducibility and the fact, that
for  $z\in D$ the probabilities $T'\delta_z$ are uniformly distributed on balls of positive radii, such $\lambda $ is unique.
The Markov operator $T$ is not strongly ergodic as topological supports of invariant ergodic probability measures $\lambda $ and $\delta_u$ ($|u| = 1$)
are not disjoint (cf [20]).
However, we can still expect  that $A_m f$ converge uniformly for some
nontrivial functions $f\in C(D_1)$.
Amongst other things our study focuses on identifying
functions $f\in C(D_1)$ such that $\lim\limits_{k\to \infty} T^k f \in C(D_1) $, where the convergence is uniform on $D_1$.
\end{exa}

\begin{defn}
Markov operators $T_i : C([0,1] \mapsto C([0,1]))$ defined as,
$$
T_i f(x) = (1-x)^{i+1} \sum_{j=0}^\infty \binom{i+j}{j}x^j f(\tfrac{j}{i+j}) \ , {\rm if } \ x \in [0,1) \ ,
$$
and \ $T_i f(1) = f(1)$, \ are called {\it Meyer-K\"{o}nig} \ or \ {\it Zeller } \ operators ($i = 1, 2,  \dots $).
\end{defn}

It was proved by I. Gavrea and M. Ivan  (see [9], [10])  that Meyer-K\"{o}nig Zeller operators have convergent iterates and
\ $\lim\limits_{m\to\infty }T_i^m f (x) = f(0) ({\bf 1} - {\bf e}_1)(x) + f(1){\bf e}_1(x)$ \ uniformly for \ $x\in [0,1]$,  \
if \ $f \in C([0,1])$.

\begin{defn}
For a fixed parameter $i \in \{ 1, 2, \dots  \}$ the linear operator defined on $C([0,1])$ by
$$
\widehat{T}_i f(x) = (i+1)(1-x)^i \sum_{j=0}^\infty \binom{i+j+1}{j}x^j \int_{\frac{j}{i+j}}^{\frac{j+1}{i+j+1}} f(t) dt  \ , {\rm if } \ x \in [0,1)
$$
and \ $\widehat{T}_i f(1) = f(1) $ \ is called a {\it generalized Kantorovich} operator.
\end{defn}

After [1] and  [2] we recall that generalized Kantorovich operators are Markovian,
$C_{{\rm inv}}([0,1]) = \{ {\rm const}\} = \{ \alpha {\bf 1} : \alpha \in \mathbb{R} \}$, and $P_{\widehat{T}_i }([0,1]) = \{ \alpha \lambda + (1- \alpha )\delta_1  : \alpha \in [0,1] \}$, where \ $\lambda$ \  is the normalized Lebesgue measure on $[0,1]$. It follows that dim$C_{{\rm inv}}([0,1]) = 1 $ and \ dim$ M_{{\rm inv}}([0,1]) = 2$. Hence $\widehat{T}_i $ is not strongly ergodic. However, we
will identify all those continuous functions $f \in C([0,1])$ whose
iterates $ \widehat{T}_i^m f \rightrightarrows$ converge uniformly. It was  proved in [1] that if such
convergence holds, then
$f(1) = \int_0^1 f(t) dt $. The question is whether (Problem 5.1 in [1])
the condition $f(1) = \int_0^1 f(t) dt $ is a sufficient one.
We solve the above posed problem and our answer is affirmative.

We shall put our solution in a general context and  address a broaden question, for general compact phase spaces $X$. In particular,
we will see that the iterates  $T^m f$ of diffusion Markov operators (from our Example 1.3) are convergent uniformly if and only if $f $ is a constant function on the boundary $\mathbb{T} = \{ z\in D_1 : | z | = 1 \}$  satisfying $f \equiv \int_{D_1} f(z) d\lambda (z)$ on $\mathbb{T}$.

In order to make our considerations rigorous we now recall after Foguel [7] necessary notions, imbedding them into a general stochastic dynamics environment.
Let $(X, \mathcal{F})$ be a measurable space, where $X$ plays the role of a phase space and $\mathcal{F}$ is a fixed $\sigma-$algebra of subsets
of $X$ (we may still keep in mind that, as before,  $X$ is a compact metric space with the Borel $\sigma-$algebra \ $\mathcal{F}$, \ completed if we are given a specific measure). By a probability transition function on $(X, \mathcal{F})$ we mean a family $\{ P(x, \cdot )  : x\in X \}$ of probability
($\sigma-$additive) measures on $(X, \mathcal{F})$ such that for each fixed $A \in \mathcal{F}$ the mapping
$X \ni x \mapsto P(x, A) \in [0,1]$ is $\mathcal{F}-$measurable. If on $(X, \mathcal{F})$ we are given a fixed positive ($\sigma-$additive), so-called reference measure $\mu$, then we usually require that $P(x,A) = 0$ for $\mu-$almost all $x$, whenever $\mu(A) = 0$.
Such a stochastic transition function allows us to define a linear operator  $T : L^{\infty}(X, \mu ) \mapsto L^{\infty}(X,\mu )$. Simply
$Th(x) = \int_X h(y) P(x, dy) = \langle h , P(x, \cdot )\rangle$, which is called
(with a slight abuse of terminolgy) Markov.
On the other hand, \  $S : L^1(X, \mu ) \mapsto L^1(X, \mu )$ \ defined by $\int_A Sf d\mu (x) = \int_X P(x, A) f(x) d\mu(x)$ is called
{\it stochastic. } It is well known (we again inherit the notation and terminology from [7]) that $S' = T$, i.e. the dual of a stochastic operator  $S$
is a Markov operator $T$. From the functional analysis perspective  stochastic and Markovian operators are positive linear contractions respectively on
$L^1(X, \mu )$ or $L^{\infty}(X, \mu )$. What is more important, \  $T{\bf 1} = {\bf 1}$ \ (in $L^{\infty}(X, \mu )$) is an invariant function
and therefore $T$ is a Markov operator in the previous sense when $L^{\infty}(X, \mu )$ is identified with $C(\triangle )$, where $\triangle$ is the compact (not necessarily metric) set of all multiplicative functionals on $L^{\infty }(X, \mu )$ endowed with the *weak topology. Let us add that every stochastic operator may be extended by $\nu S(A) = \int_X P(x,A) d\nu (x)$ to $\sigma-$additive (positive) measures $\nu $ on $(X,\mathcal{F})$.
In general, a stochastic operator does not need to posses an invariant measure (function $f\in L^1(X, \mu )$).
It is always challenging to settle if  invariant elements do exist.

\begin{defn}
A stochastic linear operator $S: L^1(X, \mu) \mapsto L^1(X, \mu )$ is called {\it doubly stochastic} if $\mu $  is a $\sigma-$additive
probability measure  and $S{\bf 1} = {\bf 1}$ (equivalently $\mu S = \mu$, or simply $\mu $ is $S-$invariant).
\end{defn}

It is well known (cf [7] or [6]) that, if $S$ is doubly stochastic, then it has the
property that  $S'(L^1(X,\mu)) \subseteq L^1(X, \mu)$. In particular, $S'$ is also doubly stochastic and we informally  write $S' = T, T' = S'' = S$.

\begin{defn}
A stochastic (linear) operator $S: L^1(X, \mu) \mapsto L^1(X, \mu )$ is called {\it kernel}  if there exists a measurable function
$$
k(\cdot , \cdot ) : (X, \mathcal{F}) \otimes (X, \mathcal{F}) \mapsto [0, \infty )
$$
such that $Sf(\cdot ) = \int_X f(x) k(x, \cdot ) d\mu (x)$.
\end{defn}

\noindent Clearly $\int_X k(\cdot , y)d\mu (y) = {\bf 1}$ as  $S'{\bf 1} = {\bf 1}$. $\nu S(A) = \int_A \int_X k(x,y) d\nu(x) d\mu (y)$ is a direct extension of a kernel stochastic operator $S$ to measures. Also notice that a $\sigma-$additive  measure
$\mu $ is $S$ invariant (i.e. $\mu S = S$) if and only if $\int_X k(x,y) d\mu(x) = 1$ for $\mu-$almost all $y\in X$.
In particular, a kernel positive operator is doubly stochastic if and only if
$$
\int_X k(\cdot , y) d\mu(y) = \int_X k(x , \cdot )d\mu (x) = {\bf 1} .
$$

Dynamic (ergodic, asymptotic) properties of stochastic (doubly stochastic) operators appear as a
leading topic in many monographs and
articles. We direct potential readers to [6], [7] or [15] as the best study materials from the perspective of
our work here. We
shall apply arguments and follow notations  presented in relatively recent articles [5], [11] and classic ones [8], [12], [16], [17], and [23].

\begin{defn}
A stochastic operator $S : L^1(X, \mu) \mapsto L^1(X, \mu )$ is called {\it asymptotically stable} if there exists
$g \in L^1(X, \mu )$ such that  for every $f \in L^1(X, \mu )$ we have
$$
\lim_{m\to \infty} \| S^m f - (\int_X f d\mu )g \|_1 = 0 \ .
$$
\end{defn}
We obviously obtain that $g$ is an $S-$invariant density ($g$ is nonnegative and $\int_X g d\mu = 1$). In particular, asymptotically stable
stochastic operators are actually doubly stochastic, with respect to a probability measure $\mu_g(A) = \int_A g d\mu $.
Without diving into details (cf [5], [7]) a kernel doubly stochastic operator is asymptotically stable if and only if
the {\it deterministic} $\sigma-$field $\Sigma_d(T)$ is trivial. Let us recall that $\Sigma_d(T) $ consists of all $ D\in \mathcal{F} $ such that
\ $T'^m{\bf 1}_D = {\bf 1}_{D_m}$ \ for all natural $m$.  When $T$ is doubly stochastic then it is a sub-$\sigma-$field of $\mathcal{F}$
(cf [7], pages 7--8  and [16]).

\section{Result}

The main result of our paper is the following:

\begin{thm}
Let $X $ be a compact Hausdorff  space and $T : C(X) \mapsto C(X)$ be a Markov linear operator satisfying:
\begin{itemize}
\item[(1)] there exists an open and dense subset $X_{\lambda } \subseteq X$ such that $X_{\lambda } \ni x \mapsto
T'\delta_x \in \mathcal{P}(X_{\lambda })$ is continuous for the total variation norm on $\mathcal{P}(X_{\lambda})$,
\item[(2)] there exists a $T-$invariant probability (Radon) measure $\lambda \in \mathcal{P}(X_{\lambda})$ such that
$\lim\limits_{m\to \infty} \| T'^m\nu - \lambda \|_{{\rm TV}} = 0$ for all probability measures $\nu $ concentrated on $X_{\lambda}$ and
$\lambda (V) > 0$ for all nonempty open $V \subseteq X_{\lambda}$,
\item[(3)] $T : C(X\setminus X_{\lambda}) \mapsto C(X\setminus X_{\lambda }) $ is a strongly ergodic Markov operator,
\item[(4)] there exists an open set $U \supseteq X\setminus X_{\lambda }$ such that $\sup_{u \in U} \| T'^2\delta_u - T'\delta_u \|_{{\rm TV}} < 2$.
\end{itemize}
Then, for a fixed function $f \in C(X)$ the iterates $T^mf$ converge uniformly on $X$ if and only if
$$
\lim_{m\to \infty} \frac{I + T + \cdots T^{m-1}}{m} f(x) = \int_X fd\lambda \ ,
$$
\noindent for all $x\in X\setminus X_{\lambda }$.
\end{thm}
Before we proceed to the proof we clarify that $\| \cdot \|_{\rm TV} $ is the total variation norm, i.e., 
$\| \nu \|_{\rm TV}
= \sup \{ |\int_X f d\nu | : \| f \| \leq 1 \}$. The Assumption (1) implies in particular, that $\| T'\delta_x - T'\delta_y \|_{\rm TV} \to 0$,
when $y\to x$ on $X_{\lambda}$. It follows from (1) and (2) that the set $X_{\lambda}$ is $T-$invariant (even though it is not closed).
The Assumption (3) implies that the set $X\setminus X_{\lambda}$ is compact and $T-$invariant.

\ \

Proof. By (2) for each $f\in C(X)$ and every $x\in  X_{\lambda}$ we have
$$
\lim_{m\to \infty} T^m f(x) = \lim_{m\to \infty} \langle f, T'^m \delta_x \rangle = \int_{X_\lambda } f d\lambda = \int_X f d\lambda \ .
$$
In other words, $T^m f (\cdot)$ converges pointwise on $X_{\lambda }$ to $\int_X d\lambda$.
Let $V$ be an open set containing $X\setminus X_{\lambda}$. Hence $X\setminus V$ is a compact subset of $X_{\lambda }$. As  the mapping
$x\mapsto T'\delta_x $ is continuous, thus it is uniformly continuous on $X\setminus V $ for the norm $\| \cdot \|_{\rm TV}$. In particular,
the set $M_{X\setminus V} = \{ T'\delta_x : x\in X \setminus V \}$
is a compact subset of $\mathcal{P}(X)$ in the norm $\| \cdot \|_{\rm TV}$. If follows from the Arzela-Ascoli theorem that
the family $\mathcal{A}_V = \{  Tf\upharpoonright_{X\setminus V} : f\in C(X) , \ \| f \|\leq 1 \}$ is relatively
compact for the sup norm $\| \cdot \|$ (indeed $Tf\upharpoonright_{X\setminus V}$ have norms bounded by 1 and are uniformly equicontinuous).
It follows that for every continuous $f\in C(X)$ the iterates $T^m f$ converge to $\int_X f d\lambda$ uniformly on $X\setminus V$, because they
converge pointwise. As $\mathcal{A}_V$ is relatively norm compact we actually have
$$
\lim_{m\to \infty} \sup_{\| f \| \leq 1} \sup_{x\in X\setminus V} |(T- I)T^mf (x)| = 0 \ .
$$
Now substitute $V = U$ and notice that
$$
\begin{aligned}
& \lim_{m\to\infty} \| T^m(T - I)\|  =
\lim_{m\to\infty} \sup_{\| f \| \leq 1} \sup_{x\in X} | T^m(T - I) f(x) | = \\
& \leq \limsup_{m\to\infty} \sup_{x\in U} \sup_{\| f \| \leq 1} |T^m (T - I) f(x) | + \limsup_{m\to\infty} \sup_{x\in X\setminus U} \sup_{\| f \| \leq 1} |T^m (T - I) f(x) | \\
& = \limsup_{m\to\infty} \sup_{x\in U} \sup_{\| f \| \leq 1} | \langle T^{m-1} f ,  T'^2\delta_x -  T'\delta_x \rangle  | \\
& \leq \limsup_{m\to\infty} \sup_{x\in U} \| T^{m-1} \| \| T'^2\delta_x  - T'\delta_x \|_{\rm TV} \leq \sup_{x\in U} \| T'^2\delta_x  - T'\delta_x \|_{\rm TV} \
< 2 \ ,
\end{aligned}
$$
\noindent by  Assumption (4). Now we shall apply the notable 0 - 2 law (see [8],  [16], [17], [18] and [23]).
For us the most practical version comes from [23] (Theorem 2.1). However, alternatively we could have turned to Theorem 3 from [17].
In any case, for every $f\in C(X)$, $\| f \| \leq 1$ and all $m\geq 1 $ we have
$$
\sup_{x\in X} |T^m (T - I) f(x)| \leq  \sup_{x\in X} |T(T - I) f(x)|   = \sup_{x\in X} |\langle f, T'^2\delta_x  - T'\delta_x \rangle | < 2 \ .
$$
It follows that the condition (iii) of Theorem 2.1 from [23] is satisfied, but therein it is equivalent to the
above
condition  (i) so we get
$$
\lim_{m\to \infty} \| T^m(T - I) f \| = 0 \ , {\rm whenever } \ f\in C(X).
$$
Now let $f\in C(X)$, $\| f \| \leq 1$ be a continuous function so that for all $x\in X \setminus X_\lambda $ we have $\frac{1}{N}\sum_{j=0}^{N-1} T^j f(x) \to \int_X f d\lambda $
(uniformly on $X\setminus X_{\lambda }$ by the Assumption (3)). It has been already proved  that \
$\lim_{j\to \infty }T^jf(x) = \int_X f d\lambda , $ \ pointwise on $X_\lambda$. It follows that $\frac{1}{N}\sum_{j=0}^{N-1} T^j f(x) \to \int_X f d\lambda =$ const for all $x\in X$. By the classical ergodic theorem (concerning  Markov operators on compact spaces, cf Proposition 1.1, Chapter 5 in [15]),  we obtain that $T$ is strongly ergodic (pointwise convergence of Cesaro means to a continuous limit  is equivalent to uniform convergence; see also [12]). Given an arbitrary $\varepsilon > 0$  we find $N_{\varepsilon}$ such that
$$
\| \frac{1}{N}\sum_{j=0}^{N-1} T^j f - (\int_X f d\lambda) {\bf 1} \| < \varepsilon \ \ {\rm for \ all} \ \ N\geq N_{\varepsilon} \ .
$$
$T$ is a linear contraction and $\lambda $ (or simply ${\bf 1}$) is $T-$invariant,
and thus for every fixed natural $m$  we have $ \| \frac{1}{N_{\varepsilon }}\sum_{j=0}^{N_{\varepsilon}-1} T^{j+m} f - (\int_X f d\lambda) {\bf 1} \| < \varepsilon $.  Besides this, for a fixed $k$ we have:
$$
\begin{aligned}
\lim_{m\to\infty} \| T^m (T^k - I) f\| &= \lim_{m\to\infty} \| T^m (T^{k-1} + T^{k-2} + \cdots + T + I)(T  - I) f\| \\
&\leq \sum_{j=0}^{k-1} \lim_{m\to\infty} \| T^{m+j} (T - I) f \| \leq \sum_{j=0}^{k-1} \lim_{m\to\infty} \| T^j \| \| T^m (T - I) f \| \\
& \leq k  \lim_{m\to\infty} \| T^m (T - I) f \| = 0 \ .
\end{aligned}
$$
Let $m_{\varepsilon} \geq N_{\varepsilon}$ be such that
$$
\| T^m (T^k - I ) f \| < \varepsilon  \ \ {\rm for \ all } \ \ k = 0, 1, \dots , N_{\vep } -1 \ ,
$$
\noindent whenever $m\geq m_{\varepsilon}$. For such $m$ we obtain
$$
\begin{aligned}
&\| T^m f - (\int_X f d\lambda) {\bf 1} \| \\
& \leq \| T^m f - \frac{1}{N_{\varepsilon }}\sum_{j=0}^{N_{\varepsilon }-1} T^{j+m} f \| + \| \frac{1}{N_{\varepsilon }}\sum_{j=0}^{N_{\varepsilon }-1} T^{j+m} f  - (\int_X f d\lambda) {\bf 1} \| \\
& < \| \frac{1}{N_{\varepsilon }}\sum_{j=0}^{N_{\varepsilon }-1} (T^m - T^{j+m} f) \| + \varepsilon \leq  \frac{1}{N_{\varepsilon }}\sum_{j=0}^{N_{\varepsilon }-1} \|T^m (I - T^{j} ) f \|  + \vep  \\
&  < \frac{1}{N_{\varepsilon}} N_{\varepsilon } \varepsilon + \vep  = 2 \vep \ .
\end{aligned}
$$
\noindent We have proved that the iterates $T^mf(\cdot )$ converge uniformly on $X$ to $(\int_X f d\lambda){\bf 1}$.

The proof in the opposite direction is self-evident. \ \ \ \ \ \ \ \ \ \ \ \ \ \ \ \ \ \ \ \ \ \ \ \ \ \ \ \ \ \ \ $\blacksquare $

\section{Affirmative answer to  Acu and Rasa conjecture }

The next theorem is a straightforward application of our Theorem 2.1. It solves Acu's and Rasa's
question (see [1] and [2]) on
asymptotic behaviour of iterates of generalized Kantorovich operators. Here the phase space is
$X = [0,1]$ with the standard modulus metric. We fix
a natural parameter $i$ and then we shall
investigate which functions $f$ have uniformly convergent iterates  $\Ki^m f$.

\begin{thm}
Let $X = [0,1]$, $i \in \mathbb{N}$ and
$$
\widehat{T}_i f(x) = (i+1)(1-x)^i \sum_{j=0}^\infty \binom{i+j+1}{j}x^j \int_{\tfrac{j}{i+j}}^{\tfrac{j+1}{i+j+1}} f(t) dt  \ , {\rm if } \ x \in [0,1) \ ,
$$
\noindent be the modified Kantorovich Markov operator \ $ \widehat{T}_i : C([0,1]) \mapsto C([0,1])$, \ where $\widehat{T}_i f(1) = f(1)$ \ ($i$ is a fixed natural parameter).
Then, 
the iterates $\widehat{T}_i^m f$ are uniformly convergent on $[0,1]$ if and only if $\int_0^1 f(u) du = f(1)$. For such $f$ we simply have
$$
\lim_{m\to\infty } \|\widehat{T}_i^m f - (\int_0^1 f(u) du){\bf 1} \| = 0 \ .
$$
\end{thm}
Proof.
It has already been mentioned (proved in [1] by Acu and Rasa) that invariant probability measures are of the form  $P_{\Ki}([0,1]) = \{ \alpha \lambda + (1-\alpha )\delta_1  : 0 \leq \alpha \leq 1 \}$. The only continuous invariant functions are constants. We find that $\Ki $ is a kernel stochastic operator
on $L^1([0,1), \lambda )$ with transition density function
$$
(\ast) \ \ \ \ \ \ \frac{\Ki'\delta_x}{d\lambda}(\cdot ) = (i+1)(1-x)^i\sum_{j=0}^{\infty }\binom{i + j + 1}{j} x^j {\bf 1}_{[\tfrac{j}{i+j} , \tfrac{j+1}{i+j+1} )}(\cdot) \ .
$$
We reiterate that $\Ki'\delta_1 = \delta_1$ and the Lebesgue probability measure $\lambda $ on $[0,1]$ is $\Ki-$invariant (clearly $\lambda ([0,1)) = 1$ but we do not claim  $\Ki$ being kernel, for the reference measure $\lambda$, on the whole $[0,1]$).
If one insists, we can introduce another invariant reference measure $\frac{1}{2}\lambda + \frac{1}{2}\delta_1$,
but then in consequence we instantly destroy the
strong Feller property. The operator $\Ki$ is doubly stochastic on $L^1([0,1), \lambda )$. Let us explicitly  write its kernel
$$
k_i(x, y) = (i+1)(1-x)^i \sum_{j=0}^{\infty } \binom{i+j+1}{j} x^j {\bf 1}_{[\tfrac{j}{i+j} , \tfrac{j+1}{i+j+1} )}(y) \geq 0 \ ,
$$
\noindent where $0 \leq x, y < 1$. It still remains to verify that $\int_0^1 k_i(x, y) dx = 1 = \int_0^1 k_i(x, y)dy$, but we skip these elementary  calculations.
When we look on $[0,1) = X_{\lambda } \ni x \mapsto \Ki'\delta_x \in L^1([0,1) , \lambda )$ and exploit the representation ($\ast $) we find that this mapping is $L^1$-norm continuous (the Assumption (1) holds). Let us comment that the Assumption (1) means that $\Ki$ is in fact strongly Feller and the image $\Ki (h) $ is a continuous function on $[0,1)$, for every measurable bounded $h$ (cf [11]), what easily follows from the representation ($\ast$). Moreover,
$\Ki'(M([0,1))) \subseteq L^1([0,1), \lambda )$. The convergence (2) holds if and only if (cf [5], Corollary 3)
the doubly stochastic operator $\Ki' : L^1([0,1), \lambda )  \mapsto L^1([0,1), \lambda ) $  has the trivial deterministic $\sigma-$field
$\Sigma_d(\Ki')$, or explicitly  $\Sigma_d(\Ki') = \Sigma_d(\Ki) = \{ \emptyset , [0,1) \}$. For this, let
a measurable set $D\subseteq [0,1)$ be such that
$0 < \lambda (D) $ and
$\Ki^m {\bf 1}_D = {\bf 1}_{D_m}$
is satisfied 
for all natural $m$. By the strong Feller property we easily find that
$$
\Ki : L^{\infty }([0,1), \lambda ) \mapsto C_b([0, 1)) \ .
$$
Therefore,
$C_b([0,1)) \ni \Ki^m {\bf 1}_D = {\bf 1}_{D_m} = {\rm const } $ as $[0,1)$ is connected. From $\lambda(D_m) = \lambda (D) > 0$ ($\Ki$ is doubly stochastic) it follows that $\lambda (D) = 1$. We have obtained
that $\Sigma_d(\Ki ) = \{  \emptyset , [0,1 ) \}$ is trivial. By [5] (Corollary 3)
the operator $\Ki'$ (also $\Ki $) is asymptotically stable on $L^1([0,1), \lambda )$. For each $g\in L^1([0,1), \lambda )$ we have
$\lim\limits_{m\to \infty} \Ki'^m g = (\int_{[0,1)} g d \lambda ){\bf 1}$ in the $L^1$ norm. Let $\nu \in \mathcal{P}([0,1))$. As $\Ki'\nu \in L^1([0,1), \lambda )$ we instantly get
$$
\lim_{m\to \infty } \| \Ki'^m \nu - \lambda \|_{\rm TV} = \lim_{m\to \infty } \| \Ki'^{m-1}(\tfrac{\Ki'\nu}{d\lambda } ) - {\bf 1} \|_{L^1} = 0 \ .
$$
Here $[0,1] \setminus X_{\lambda } = [0,1] \setminus [0,1) = \{ 1 \}$. As $\Ki'\delta_1 = \delta_1$ it follows that $\Ki $ is the identity
operator on $C(X\setminus X_{\lambda})$ and the Assumption (3) is trivially satisfied. One can be tempted to finish
by applying  Jamison's Theorem from [12]. However, it is not applicable as there the operator is assumed to be irreducible on the whole phase space $X$. Our operator $\Ki$ has two ergodic invariant probabilities, so it is impossible
to comply with this assumption.

We must verify whether the Assumption (4) holds. This is the most technical and laborious part of our work.
On the other hand, it merely hinges on elementary manipulations. We shall estimate
$\| \Ki'^2\delta_x - \Ki'\delta_x \|_{ \rm TV}$ for $x$ close to $1$. Our target is to show that
$$
\inf \{ \| \Ki'^2 \delta_x \wedge \Ki'\delta_x  \|_{\rm TV} : x > 1 - \widetilde{\varepsilon }  \}  > \varepsilon
$$
\noindent for some $\vep , \tvep > 0 $, where $\wedge$ denotes the standard lattice minimum in $M([0,1))$.
The detailed description of the functions $\alpha_j(x)$ and $\beta_j(x)$, introduced below,
are not really necessary. Nevertheless, they help us understanding the behaviour of $\Ki'\delta_x$ on
component segments $[\frac{j}{j+i} , \frac{j+1}{i+j+1} )$, when $x\to 1^-$.
So, let us start with ($\ast $) and introduce
$$
\alpha_j(x) = (i+1)\binom{i + j + 1}{j}(1-x)^i x^j  \ ,
$$
\noindent  $x\in [0,1)$. Notice that $\alpha_0(x) = (i+1)(1-x)^i$ for all $x\in [0, 1).$
$$
\begin{aligned}
\Ki'\delta_x &= (i+1)(1-x)^i\sum_{j=0}^{\infty }\binom{i + j + 1}{j} x^j {\bf 1}_{[\tfrac{j}{i+j} , \tfrac{j+1}{i+j+1} )}(\cdot)
= \sum_{j=0}^\infty \alpha_j(x) {\bf 1}_{[\tfrac{j}{i+j} , \tfrac{j+1}{i+j+1} )}(\cdot) \ .
\end{aligned}
$$
We easily get
$$
\begin{aligned}
\frac{\alpha_{j+1}(x)}{\alpha_j(x)} &= \frac{(i+1)\binom{i + j + 2}{j+1}(1-x)^i x^{j+1}}{(i+1)\binom{i + j + 1}{j}(1-x)^i x^j}
= (1 + \tfrac{i+1}{j+1})x \ (\to x \ {\rm when } \ j\to \infty ) .
\end{aligned}
$$
In particular, $\alpha_j(x) =  \prod\limits_{k=1}^j (1 + \tfrac{i+1}{k}) x^j\alpha_0(x)  = (i+1)(1-x)^ix^j\prod\limits_{k=1}^j (1 + \tfrac{i+1}{k}) .$ Moreover, $\frac{\alpha_{j+1}(x)}{\alpha_j(x)} \geq 1 \Leftrightarrow j \leq \frac{(i+2)x - 1}{1-x}$. Because we focus on $x$ close to $1$, thus  without
loss of generality we may consider $x > \frac{1}{i+2}$ exclusively. Thus for a fixed $\frac{1}{i+2} < x < 1$ the sequence
$\alpha_j(x)$ increases for $j = 0, 1, \dots ,\lfloor\frac{(i+2)x - 1}{1-x}\rfloor , \lfloor \frac{(i+1)x }{1-x}\rfloor $ and decreases for $j\geq  \lfloor \frac{(i+1)x }{1-x} \rfloor $.
We easily find $\lim\limits_{j\to \infty} \alpha_j(x) = 0$ (asymptotically with the geometric rate depending on $x$). Its maximum is obtained
for $j_{\rm max} =  \lfloor\frac{(i+2)x - 1}{1-x}\rfloor + 1 = \lfloor \frac{(i+1)x}{1-x} \rfloor $.
Now
$$
\begin{aligned}
\alpha_{j_{\rm max}} (x) &= (i+1)\frac{(i+j_{\rm max} + 1)!}{(i+1)!j_{\rm max}!} (1-x)^i x^{j_{\rm max}} \\
&= \frac{(1-x)^i}{i!}j_{\rm max}^{i+1}(1 + \tfrac{1}{j_{\rm max}})\cdots (1 + \tfrac{i+1}{j_{\rm max}})x^{j_{\rm max}} \\
&\geq \frac{((i+2)x - 1)^{i+1}}{(1-x) i!} x^{\tfrac{(i+1)x}{1-x}} \ .\\
\end{aligned}
$$
We estimate for  $x \to 1^-$
$$
x^{\tfrac{(i+1)x}{1-x}} = \left( (1 - (1-x))^{\tfrac{1}{1-x}}  \right) ^{(i+1)x} \to \left( \frac{1}{\rm e}\right)^{i+1} .
$$
Denote $\tvep > 0$ such that  $\alpha_{j_{\rm max}} (x) \geq \frac{1}{2} \frac{(i+1)^i}{i!} \left( \frac{1}{\rm e} \right)^{i+1} \frac{1}{1-x}$ when $x > 1 - \tvep$.

One more component function of $\Ki'\delta_x $  is
$$
\beta_j (x) = \Ki'\delta_x ([\frac{j}{i+j} , \frac{j+1}{i + j + 1})) \ .
$$
\noindent  After elementary calculations
$$
\beta_j (x) = \alpha_j(x) \left( \tfrac{j+1}{i + j +1} - \tfrac{j}{i+j} \right) \\
= \binom{i + j -1}{j} (1-x)^i x^j \ .
$$
We generally  have  $\beta_0(x) = (1-x)^i$ and in the particular case when $i=1$ we have $\beta_j(x) = x^j (1-x)$.
Now we describe the dynamics of $\beta_j(x)$ for general $i\ge 1$. The ratio
$$
\frac{\beta_{j+1}(x)}{\beta_j(x)} = \frac{\binom{i+j}{j+1}(1-x)^ix^{j+1}}{\binom{i+j-1}{j} (1-x)^ix^j} \\
= \left(  1 +\tfrac{i-1}{j+1}  \right)x \ \ {\rm for } \ \ j = 0, 1, \dots
$$
Hence,
$$
\beta_j(x) = \prod\limits_{k=1}^j (1 + \tfrac{i-1}{k})\cdot x^j \beta_0(x) =  \prod\limits_{k=1}^j (1 + \tfrac{i-1}{k})\cdot x^j (1-x)^i \ ,
$$
\noindent here we keep the standard convention $\prod\limits_{k=1}^0 \dots = 1 $.
Now $\frac{\beta_{j+1}(x)}{\beta_j(x)} \geq 1 $ if and only if $j+1 \leq \frac{(i-1)x}{1-x}$. If $i =1$ we simply  have
$\frac{\beta_{j+1}(x)}{\beta_j(x)} = x$.  It follows that for $i = 1$ the sequence $\beta_j(x)$ is strictly decreasing and $\beta_0(x)
= (1 -x) > x^j (1-x)$ (if $i=1$ then  argmax$_j \beta_j(x) = 0$). Generally we have
$$
{\rm argmax}_j \beta_j (x) = \left\lfloor \frac{(i-1)x}{1-x}\right\rfloor \ .
$$
We shall not exploit this fact but it is worth 
realizing that, if $i\geq 2$,  then  we have
$$
{\rm argmax}_j \alpha_j(x) - {\rm argmax}_j \beta_j(x) \leq \frac{(i+1)x}{1-x} - \frac{(i-1)x}{1-x} + 1 = \frac{1 + x}{1-x}
$$
and
$$
{\rm argmax}_j \alpha_j(x) - {\rm argmax}_j \beta_j(x) \geq  \left(\frac{(i+1)x}{1-x} -1\right) - \frac{(i-1)x}{1-x}  = \frac{3x-1}{1-x} \ .
$$
The discrepancy $\frac{1 + x}{1-x} - \frac{3x-1}{1-x} = 2$ is not considerably high (except for $i=1$, when these features completely fritter away).
We notice that $\frac{(i+1)x}{1-x} \to \infty $ and $\frac{(i-1)x}{1-x} \to \infty $ when $x\to 1^-$ (excluding the trivial case $i=1$).

However, for simplicity (and to cover the case $i=1$) it is better to take another index.  Namely let
$j_x$ be the unique nonnegative integer  $j \in \{ 0, 1, \dots \}$ such that $x\in [\tfrac{j}{i+j} , \tfrac{j+1}{i+j+1})$ (we can find
that $\frac{ix}{1-x} -1 < j_x \leq \frac{ix}{1-x}$). Clearly $[0, 1) \ni x \mapsto j_x \in \{ 0, 1, \dots \}$ is nondecreasing.
We have
$$
{\rm argmax}_j\beta_j(x) = \left\lfloor \frac{(i-1)x}{1-x} \right\rfloor < j_x \leq \frac{ix}{1-x} \leq \left\lfloor
\frac{ix + x}{1-x} \right\rfloor \leq  {\rm argmax}_j\alpha_j(x) \ .
$$
We find that $j_x $ differs from argmax$_j\alpha_j(x)$ and argmax$_j\beta_j(x)$.
Obviously $j_x \to \infty $ if and only if $x\to 1^- $.

We remember (Calculus I) that $(1+\frac{1}{t})^t $ , $(1 - \frac{1}{t})^t$ converge
monotonically (the latter from $t \ge 1$) respectively to ${\rm e}$ or $\frac{1}{\rm e}$
when $t$ increases to $\infty$.

Let us fix $\vep > 0 $ (small) and choose $r = r_{\vep} \in (\tfrac{\vep}{1 + \vep } , 1)$. We easily find that $(1+\vep )(1-r) < 1$.  Define $u_{\vep}$ so close to $1^-$ that
$(1 + \tfrac{1}{t})^t > {\rm e}^{1-\vep}$ and $(1 - \tfrac{1}{t})^t > {\rm e}^{-1-\vep } = (1-\kappa_{\vep})\tfrac{1}{\rm e}$, whenever
$t \geq \frac{\lfloor(1-r)j_x\rfloor}{i}$ for all $x\in (u_{\vep}, 1)$. In the sequel we abbreviate $\kappa = \kappa_{\vep}$ when it appears in
formulae.
We set
$$
C_x = \sum_{j = \lfloor(1-r) j_x\rfloor}^{j_x} \beta_j(x) > 0 \ .
$$
We work toward showing that $C_x$ are uniformly separated from $0$ for $x$ being close enough to $1^-$ (how close is sufficient we shall learn
in due course).
We have
$$
C_x = \sum_{j = \lfloor(1-r) j_x\rfloor}^{j_x} \left( \prod_{k=1}^j (1 + \tfrac{i-1}{k})\right) (1-x)^i x^j \ .
$$
It is well known that $\prod_{k=1}^j (1 + \tfrac{i-1}{k}) \sim j^{i-1}$ when $j\to \infty $. Let $\upsilon_{\kappa} < 1$ be so close to $1^-$
that $\prod_{k=1}^j (1 + \tfrac{i-1}{k}) > (1-\kappa )j^{i-1}$, whenever $j \geq \lfloor (1-r)j_x\rfloor $ for all $\upsilon_{\kappa} < x < 1 .$
\ Now for all $1 > x > \max \{  u_{\vep} ,  \upsilon_{\kappa} \}$ we have
$$
\begin{aligned}
C_x &\geq \sum_{j = \lfloor(1-r) j_x\rfloor}^{j_x} (1-\kappa )j^{i-1}x^j(1-x)^i = (1-\kappa)(1-x) \sum_{j = \lfloor(1-r) j_x\rfloor}^{j_x} x^j(j(1-x))^{i-1} \\
&\geq (1-\kappa)(1-x) \sum_{j = \lfloor(1-r) j_x\rfloor}^{j_x} x^j \left( j\cdot \tfrac{i}{i +j_x +1} \right)^{i-1}.
\end{aligned}
$$
We estimate (for \  $j\in \{\lfloor(1-r) j_x\rfloor , \lfloor(1-r) j_x\rfloor +1, \dots , j_x \}$)
$$
\begin{aligned}
&\left( j\cdot \tfrac{i}{i +j_x +1} \right)^{i-1}  = \left( \frac{j}{j_x}\right)^{i-1}\cdot \left(\tfrac{i}{1 + \tfrac{i + 1}{j_x} } \right)^{i-1}
\geq  \left(  \tfrac{(1-r)j_x -1}{j_x} \right)^{i-1}\cdot \tfrac {i^{i-1}}{\left(  1 + \tfrac{i+1}{j_x}\right)^{i-1}} \\
&= \left(  1-r -\tfrac{1}{j_x} \right)^{i-1}\cdot \frac{i^{i-1}}{\left(  1 + \tfrac{i+1}{j_x}\right)^{i-1}}
\geq (1-\widetilde{r})^{i-1}i^{i-1} \
\end{aligned}
$$
\noindent for some $\widetilde{r} \in (r, 1)$ ($x\ > \Psi$ described below). Indeed, it is sufficient to take $x$ close enough to $1^-$, so $j_x$ and $\widetilde{r} $
are linked to each other through $j_x > \frac{1}{\widetilde{r} - r}$ and $\left(  1 + \frac{i+1}{j_x} \right)^{i-1}$ may be assumed to be close to $1$.
Notice that $\lim\limits_{x\to 1^-}\left(  1 + \frac{i+1}{j_x} \right)^{i-1} = 1.$ Hence, for some $\Psi < 1$ the state $x\in ( \Psi , 1)$  meets the required inequalities.
Now
$$
\begin{aligned}
C_x &\geq (1-\kappa)(1-x) \sum_{j = \lfloor(1-r) j_x\rfloor}^{j_x} x^j \left( j\cdot \tfrac{i+1}{i +j_x +1} \right)^{i-1} \\
&\geq (1-\kappa)(1-\widetilde{r}) i^{i-1} (1-x) x^{\lfloor (1-r)j_x \rfloor} \left( 1 + x + \cdots + x^{j_x - \lfloor (1-r)j_x\rfloor } \right) \\
&= (1-\kappa )(1-\widetilde{r} )i^{i-1} (x^{\lfloor (1-r)j_x\rfloor } - x^{j_x+1}) \ .
\end{aligned}
$$
It remains to estimate
$$
\begin{aligned}
x^{\lfloor (1-r)j_x\rfloor } &- x^{j_x+1} \geq \left( \tfrac{j_x}{i+j_x}  \right)^{\lfloor (1-r)j_x \rfloor} -
\left( \tfrac{j_x +1}{i+j_x+1}  \right)^{j_x +1} \\
&\geq \left( 1 - \tfrac{i}{i+j_x} \right)^{(1-r)j_x} - \left( 1 - \tfrac{i}{i+j_x+1} \right)^{j_x +1} \\
&= \left[ \left( 1 - \tfrac{1}{1+\tfrac{j_x}{i}} \right)^{1+\tfrac{j_x}{i}}  \right]^{\tfrac{(1-r)j_x}{1+\tfrac{j_x}{i}}} -
\left[ \left( 1 - \tfrac{1}{1+\tfrac{j_x +1}{i} } \right)^{\tfrac{i +j_x +1}{i}}   \right]^{\tfrac{(j_x + 1)i}{i +j_x +1}} \\
&\geq \left(  {\rm e}^{-1 - \vep}\right)^{\tfrac{1-r}{\tfrac{1}{i} + \tfrac{1}{j_x}}} -   {\rm e}^{-i\tfrac{j_x + 1}{i + j_x +1 }}  \\
&\geq {\rm e}^{-(1+\vep)(1-r)i} - {\rm e}^{-i\eta}  > 0
\end{aligned}
$$
as $(1+\vep )(1-r) < (1+\vep)\tfrac{1}{1+\vep} = 1 $ and $\tfrac{j_x +1}{i+ j_x + 1} > \eta > (1+\vep)(1-r) $ when
$ 1 > x > \Theta $, where $\Theta \in (0,1)$ is appropriately large.  Completing this part of our proof we have
$$
C_x \geq C  = (1-\kappa)(1-\widetilde{r})^{i-1}i^{i-1}\left( {\rm e}^{-i (1+\vep)(1-r)} - {\rm e}^{-i \eta} \right) > 0 ,
$$
\noindent for all $x\in \left[\max\{  u_{\vep} , \upsilon_{\vep} , w, \Theta ,  \Psi \}   , 1 \right), $ where $C$ does not depend on such $x$.

Now define $\gamma_l(x) = \Ki'^2\delta_x([ \tfrac{l}{i+l }, \tfrac{l+1}{i + l + 1}))$. Let us recall that
for each probability measure $\nu $ on $[0,1)$ we have a representation (in fact it is a density $\frac{\Ki'\delta_x}{d\lambda}$)
$$
\begin{aligned}
\Ki'\nu &= (i+1)\sum_{l=0}^{\infty } \binom{i+l+1}{l}\int_0^1 (1-y)^iy^l d\nu(y) {\bf 1}_{[\tfrac{l}{i+l}, \tfrac{l+1}{i + l + 1})}(\cdot )  .
\end{aligned}
$$
 We evaluate (thinking on $\nu = \Ki'\delta_x$)
$$
\begin{aligned}
&\int_0^1 (1-y)^iy^l d\Ki'\delta_x (y) =\sum_{j=0}^{\infty }(i+1)\binom{i+j+1}{j} (1-x)^i x^j \int_{\tfrac{j}{i+j}}^{\tfrac{j+1}{i+j+1}} (1-y)^iy^l dy\\
&\geq \sum_{j=0}^{\infty} (i+1)\binom{i+j+1}{j} (1-x)^ix^j \left(\frac{j}{i+j}\right)^l\left(\tfrac{i}{i+j+1}\right)^i \frac{i}{(i+j)(i+j+1)} \\
&= \sum_{j=0}^{\infty} \left( \prod_{k=1}^j (1 + \tfrac{i-1}{k}) \right)(1-x)^i x^j \left( \tfrac{j}{i+j}\right)^l\left( \tfrac{i}{i+j+1}\right)^i \ .
\end{aligned}
$$
Hence
$$
\begin{aligned}
&\gamma_l(x) = \Ki'^2\delta_x([ \tfrac{l}{i+l }, \tfrac{l+1}{i + l + 1})) \\
&\geq (i+1)\binom{i+l+1}{l} \sum_{j=0}^\infty \left( \prod_{k=1}^j (1 + \tfrac{i-1}{k})\right)(1-x)^i x^j \left( \tfrac{j}{i+j}\right)^l\left( \tfrac{i}{i+j+1}\right)^i \cdot \frac{i}{(i+l)(i+l+1)} \\
&= \sum_{j=0}^\infty \left( \prod_{k=1}^j (1 + \frac{i-1}{k})\right)(1-x)^i x^j \left( \prod_{m=1}^l (1+\tfrac{i-1}{m})\right)\left( \tfrac{j}{i+j}\right)^l\left( \tfrac{i}{i+j+1}\right)^i  \ .
\end{aligned}
$$
We shall compare $\beta_l(x)$ with $\gamma_l(x)$, where $\lfloor (1-r)j_x \rfloor \leq l  \leq j_x$.  For this we shall
use the inequality $\frac{i}{i+ j_x +1} \leq 1 - x$ and consequently, via the estimation (for $j\leq j_x$)
$$
\begin{aligned}
\left( \frac{\tfrac{i}{i+j+1}}{1-x}\right)^i &\geq \left(\frac{\tfrac{i}{i + j + 1}}{1 - \tfrac{j_x }{i + j_x +1 }} \right)^i \geq\left( \frac{i}{i+1} \right)^i \ ,
\end{aligned}
$$
we approach the end considering the ratio
$$
\begin{aligned}
\frac{\gamma_l(x)}{\beta_l(x)} &\geq\frac{\sum_{j=0}^\infty \left( \prod_{k=1}^j (1 + \tfrac{i-1}{k})\right)(1-x)^i x^j \left( \prod_{m=1}^l (1+\tfrac{i-1}{m})\right)\left( \tfrac{j}{i+j}\right)^l\left( \tfrac{i}{i+j+1}\right)^i}{\left(\prod\limits_{m=1}^l (1 + \tfrac{i-1}{m})\right)\cdot  (1-x)^i\cdot x^l} \\
&=\sum_{j=0}^\infty \left( \prod_{k=1}^j (1 + \tfrac{i-1}{k})\right)x^j \left( \frac{\tfrac{j}{i+j}}{x}\right)^l \left( \frac{\tfrac{i}{i+j+1}}{1-x}\right)^i \\
&\geq \sum_{j=\lfloor (1-r)j_x \rfloor }^{j_x} \left( \prod_{k=1}^j (1 + \tfrac{i-1}{k})\right)x^j \frac{\left( \tfrac{j}{i+j} \right)^l}{ 1}\cdot \left( \frac{i}{i+1}\right)^i\\
&\geq  \left( \frac{i}{i+1}\right)^i\sum_{j=\lfloor (1-r)j_x \rfloor }^{j_x} \left( \prod_{k=1}^j (1 + \tfrac{i-1}{k})\right)x^j \left[ \left( 1 - \tfrac{i}{i+j} \right)^{i+j}\right]^{\tfrac{l}{i+j}} \\
&\geq \left( \frac{i}{i+1}\right)^i\sum_{j=\lfloor (1-r)j_x \rfloor }^{j_x} \left( \prod_{k=1}^j (1 + \tfrac{i-1}{k})\right)x^j \left[ \left( 1 - \kappa \right)\tfrac{1}{{\rm e}^i} \right]^{\tfrac{j_x}{i+(1-r)j_x}}\\
&\geq \left( \frac{i}{i+1}\right)^i\sum_{j=\lfloor (1-r)j_x \rfloor }^{j_x} \left( \prod_{k=1}^j (1 + \tfrac{i-1}{k})\right) x^j \left[ \left( 1 - \kappa \right)\tfrac{1}{{\rm e}^i} \right]^{\tfrac{1}{1-r}}\\
&= \left( \frac{i}{i+1}\right)^i \sum_{j=\lfloor (1-r)j_x \rfloor }^{j_x} \left( \prod_{k=1}^j (1 + \tfrac{i-1}{k})\right) x^j (1-\kappa)^{\tfrac{1}{1-r}}{\rm e}^{-\tfrac{i}{1-\kappa }} \\
&=  \left( \frac{i}{i+1}\right)^i (1-\kappa)^{\tfrac{1}{1-r}}{\rm e}^{-\tfrac{i}{1-\kappa }}\sum_{j=\lfloor (1-r)j_x \rfloor }^{j_x} \left( \prod_{k=1}^j (1 + \tfrac{i-1}{k})\right) x^j  \\
&\geq C \left( \frac{i}{i+1}\right)^i (1-\kappa)^{\tfrac{1}{1-r}}{\rm e}^{-\tfrac{i}{1-\kappa }} > 0 \ .
\end{aligned}
$$
Thus,
$$
\gamma_j(x) \geq  C \left( \frac{i}{i+1}\right)^i (1-\kappa)^{\tfrac{1}{1-r}}{\rm e}^{-\tfrac{i}{1-\kappa }}\beta_j(x) \ ,
$$
\noindent for all $\lfloor (1-r) j_x\rfloor \leq j \leq j_x$.
Now,
$$
\begin{aligned}
\| \Ki'^2\delta_x \wedge \Ki'\delta_x \|_{\rm TV} &\geq \sum_{j=\lfloor (1-r)j_x \rfloor }^{j_x} \gamma_j(x)\wedge \beta_j(x) \\
&\geq \left[\left( C \left( \frac{i}{i+1}\right)^i (1-\kappa)^{\tfrac{1}{1-r}}{\rm e}^{-\tfrac{i}{1-\kappa }} \right) \wedge 1 \right]\sum_{j=\lfloor (1-r)j_x \rfloor }^{j_x} \beta_j(x) \\
&\geq \left[ \left( C \left( \frac{i}{i+1}\right)^i (1-\kappa)^{\tfrac{1}{1-r}}{\rm e}^{-\tfrac{i}{1-\kappa }} \right) \wedge 1 \right] C > 0 \ ,
\end{aligned}
$$
\noindent for all $x\in (\max \{u_{\vep}, \upsilon_{\vep}, w, \Theta , \Psi \} , 1 )$.
It follows that
$$
\| \Ki'^2\delta_x - \Ki'\delta_x \|_{\rm } = 2 (1 - \|\Ki'^2\delta_x \wedge \Ki'\delta_x \|_{\rm TV}) < 2
$$
\noindent for all $x\in U$, where the set $U = (\max \{u_{\vep}, \upsilon_{\vep}, w, \Theta, \Psi  \}, 1]$ is a required open
set covering $\{ 1 \},$ satisfying 
Assumption (4). \ \ \ \ \ \ \ \ \ \ \ \ $ \blacksquare $

\ \

\end{document}